\documentclass[conference]{IEEEtran}
\usepackage[latin9]{inputenc}
\usepackage{listings}
\usepackage{units}
\usepackage{amsthm}
\usepackage{amsmath}
\usepackage{amssymb}
\usepackage{graphicx}
\usepackage{esint}

\makeatletter

\providecommand{\tabularnewline}{\\}
\theoremstyle{plain}
\newtheorem{thm}{\protect\theoremname}
\theoremstyle{definition}
\newtheorem{defn}[thm]{\protect\definitionname}
\theoremstyle{plain}
\newtheorem{conjecture}[thm]{\protect\conjecturename}
\theoremstyle{plain}
\newtheorem{lem}[thm]{\protect\lemmaname}
\theoremstyle{plain}
\newtheorem{cor}[thm]{\protect\corollaryname}
\theoremstyle{plain}
\newtheorem{prop}[thm]{\protect\propositionname}
\theoremstyle{definition}
\newtheorem{example}[thm]{\protect\examplename}

\usepackage{amsfonts}
\usepackage{cite}
\usepackage{pgf,tikz}
\usetikzlibrary{arrows}
\usepackage{xcolor}
\usepackage{units}

\makeatother

\providecommand{\conjecturename}{Conjecture}
\providecommand{\corollaryname}{Corollary}
\providecommand{\definitionname}{Definition}
\providecommand{\examplename}{Example}
\providecommand{\lemmaname}{Lemma}
\providecommand{\propositionname}{Proposition}
\providecommand{\theoremname}{Theorem}

\begin{document}

\title{Mutual information of Contingency Tables and Related Inequalities}

\author{\authorblockN{Peter Harremo{\"e}s} \authorblockA{Copenhagen Business
College\\
 Copenhagen, Denmark\\
 Email: harremoes@ieee.org} }
\maketitle
\begin{abstract}
For testing independence it is very popular to use either the $\chi^{2}$-statistic
or $G^{2}$-statistics (mutual information). Asymptotically both are
$\chi^{2}$-distributed so an obvious question is which of the two
statistics that has a distribution that is closest to the $\chi^{2}$-distribution.
Surprisingly the distribution of mutual information is much better
approximated by a $\chi^{2}$-distribution than the $\chi^{2}$-statistic.
For technical reasons we shall focus on the simplest case with one
degree of freedom. We introduce the signed log-likelihood and demonstrate
that its distribution function can be related to the distribution
function of a standard Gaussian by inequalities. For the hypergeometric
distribution we formulate a general conjecture about how close the
signed log-likelihood is to a standard Gaussian, and this conjecture
gives much more accurate estimates of the tail probabilities of this
type of distribution than previously published results. The conjecture
has been proved numerically in all cases relevant for testing independence
and further evidence of its validity is given.
\end{abstract}

\section{Choice of statistic}

We consider the problem of testing independence in a discrete setting.
Here we shall follow the classic approach to this problem as developed
by Pearson, Neuman and Fisher. The question is whether a sample with
observation counts $\left(X_{ij}\right)$ has been generated by the
distribution $Q=\left(q_{ij}\right)$ where $q_{ij}=s_{i}\cdot t_{j}$
for some probability vectors $\left(s_{i}\right)$ and $\left(t_{j}\right)$
reflecting independence. We introduce the marginal counts $R_{i}=\sum_{j}X_{ij}$
and $S_{j}=\sum_{i}X_{ij}$ and the sample size $N=\sum_{ij}X_{ij}$.
The maximum likelihood estimates of probability vectors $\left(s_{i}\right)$
and $\left(t_{j}\right)$ are given by $s_{i}=\frac{R_{i}}{N}$ and
$t_{j}=\frac{S_{j}}{N}$ leading to $q_{ij}=\frac{R_{i}\cdot S_{j}}{N^{2}}$.
We also introduce the empirical distribution $\hat{P}=\left(\frac{X_{ij}}{N}\right)$
.

Often one uses one of the Csisz{\'a}r \cite{Csiszar1963} $f$-divergences
\begin{equation}
D_{f}\left(\hat{P},Q\right)={\textstyle \sum\limits _{i,j}}q_{ij}f\left(\frac{\hat{p}_{ij}}{q_{ij}}\right).\label{Csi}
\end{equation}
 The null hypothesis is accepted if the test statistic $D_{f}\left(\hat{P},Q\right)$
is small and rejected if $D_{f}\left(\hat{P},Q\right)$ is large.
Whether $D_{f}\left(\hat{P},Q\right)$ is considered to be small or
large depends on the significance level \cite{Lehman2005}. The most
important cases are obtained for the convex functions $f(t)=n(t-1)^{2}$
and $f(t)=2nt\ln t$ leading to the \emph{Pearson }$\chi^{2}$\emph{-statistic}
\begin{equation}
\chi^{2}={\textstyle \sum\limits _{i,j}}\frac{(X_{ij}-Nq_{ij})^{2}}{Nq_{nj}}\label{Pears}
\end{equation}
and the \emph{likelihood ratio statistic }
\begin{equation}
G^{2}=2{\textstyle \sum\limits _{i,j}}X_{ij}\ln\frac{X_{ij}}{Nq_{ij}}.\label{lik}
\end{equation}
We note that $G^{2}$ equals the mutual information between the indices
$i$ and $j$ times $2n$ when the empirical distribution $\hat{P}$
is used as joint distribution over $i$ and $j$. See \cite{Csiszar2004}
for a short introduction to contingency tables and further references
on the subject.

One way of choosing between various statistics is by computing their
asymptotic efficiency \cite{Harremoes2012,Harremoes2008} but for
finite sample size such results are difficult to use. Therefore we
will turn our attention to another property that is of importance
in choosing a statistic.

\begin{figure}
\centering{}\definecolor{ffqqqq}{rgb}{1,0,0} 
\definecolor{qqqqff}{rgb}{0,0,1} 
\begin{tikzpicture}[line cap=round,line join=round,>=triangle 45,x=0.028001020538333976cm,y=0.028001020538333973cm] 
\draw[->,color=black] (-8,0) -- (148,0); 
\foreach \x in {20,40,60,80,100,120,140} 
\draw[shift={(\x,0)},color=black] (0pt,2pt) -- (0pt,-2pt) node[below] {\footnotesize $\x$}; 
\draw[->,color=black] (0,-8) -- (0,148); 
\foreach \y in {20,40,60,80,100,120,140} \draw[shift={(0,\y)},color=black] (2pt,0pt) -- (-2pt,0pt) node[left] {\footnotesize $\y$}; 
\draw[color=black] (0pt,-10pt) node[right] {\footnotesize $0$}; 
\clip(-8,-8) rectangle (148,148); 
\draw [color=qqqqff] (119,127)-- (119,111); 
\draw [color=qqqqff] (105,99)-- (105,111); 
\draw [color=qqqqff] (93,99)-- (93,88); 
\draw [color=qqqqff] (83,88)-- (83,78); 
\draw [color=qqqqff] (74,78)-- (74,69); 
\draw [color=qqqqff] (65,69)-- (65,61); 
\draw [color=qqqqff] (58,61)-- (58,54); 
\draw [color=qqqqff] (51,54)-- (51,47); 
\draw [color=qqqqff] (44,47)-- (44,41); 
\draw [color=qqqqff] (39,41)-- (39,36); 
\draw [color=qqqqff] (33,36)-- (33,31); 
\draw [color=qqqqff] (28,31)-- (28,26); 
\draw [color=qqqqff] (24,26)-- (24,22); 
\draw [color=qqqqff] (20,22)-- (20,18); 
\draw [color=qqqqff] (16,18)-- (16,15); 
\draw [color=qqqqff] (13,15)-- (13,12); 
\draw [color=qqqqff] (10,12)-- (10,9); 
\draw [color=qqqqff] (8,9)-- (8,7); 
\draw [color=qqqqff] (6,7)-- (6,5); 
\draw [color=qqqqff] (4,5)-- (4,3); 
\draw [color=qqqqff] (3,3)-- (3,2); 
\draw [color=qqqqff] (1,2)-- (1,1); 
\draw [color=qqqqff] (1,1)-- (1,0); 
\draw [color=qqqqff] (0,0)-- (0,0); 
\draw [color=qqqqff] (139,127)-- (139,143); 
\draw[color=ffqqqq,smooth,samples=100,domain=0.0:148.0] plot(\x,{(\x)}); \end{tikzpicture}\caption{\label{Fig1}Quantiles of twice mutual information vs. a $\chi^{2}$-distribution.}
\end{figure}
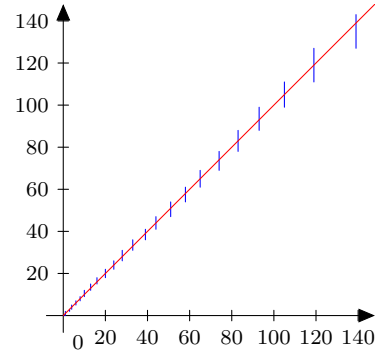
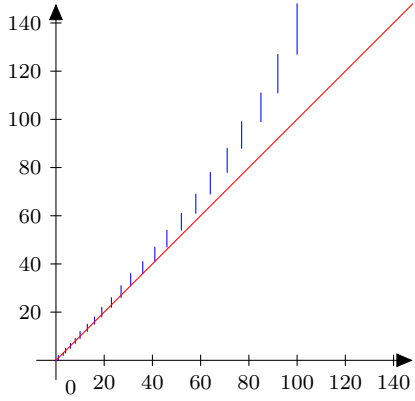
\begin{figure}
\centering{}\definecolor{ffqqqq}{rgb}{1,0,0} 
\definecolor{qqqqff}{rgb}{0,0,1} 
\begin{tikzpicture}[line cap=round,line join=round,>=triangle 45,x=0.032051282051282055cm,y=0.03205128205128205cm] 
\draw[->,color=black] (-8,0) -- (148,0); 
\foreach \x in {20,40,60,80,100,120,140} 
\draw[shift={(\x,0)},color=black] (0pt,2pt) -- (0pt,-2pt) node[below] {\footnotesize $\x$}; 
\draw[->,color=black] (0,-8) -- (0,148); 
\foreach \y in {20,40,60,80,100,120,140} 
\draw[shift={(0,\y)},color=black] (2pt,0pt) -- (-2pt,0pt) node[left] {\footnotesize $\y$}; \draw[color=black] (0pt,-10pt) node[right] {\footnotesize $0$}; 
\clip(-8,-8) rectangle (148,148); 
\draw [color=qqqqff] (92,127)-- (92,111); 
\draw [color=qqqqff] (85,99)-- (85,111); 
\draw [color=qqqqff] (77,99)-- (77,88); 
\draw [color=qqqqff] (71,88)-- (71,78); 
\draw [color=qqqqff] (64,78)-- (64,69); 
\draw [color=qqqqff] (58,69)-- (58,61); 
\draw [color=qqqqff] (52,61)-- (52,54); 
\draw [color=qqqqff] (46,54)-- (46,47); 
\draw [color=qqqqff] (41,47)-- (41,41); 
\draw [color=qqqqff] (36,41)-- (36,36); 
\draw [color=qqqqff] (31,36)-- (31,31); 
\draw [color=qqqqff] (27,31)-- (27,26); 
\draw [color=qqqqff] (23,26)-- (23,22); 
\draw [color=qqqqff] (19,22)-- (19,18); 
\draw [color=qqqqff] (16,18)-- (16,15); 
\draw [color=qqqqff] (13,15)-- (13,12); 
\draw [color=qqqqff] (10,12)-- (10,9); 
\draw [color=qqqqff] (8,9)-- (8,7); 
\draw [color=qqqqff] (6,7)-- (6,5); 
\draw [color=qqqqff] (4,5)-- (4,3); 
\draw [color=qqqqff] (3,3)-- (3,2); 
\draw [color=qqqqff] (1,2)-- (1,1); 
\draw [color=qqqqff] (1,1)-- (1,0); 
\draw [color=qqqqff] (0,0)-- (0,0); 
\draw [color=qqqqff] (100,127) -- (100,148); 
\draw[color=ffqqqq,smooth,samples=100,domain=0.0:147.99999999999997] plot(\x,{(\x)}); \end{tikzpicture}\caption{Quantiles of the $\chi^{2}$-statistic vs. a $\chi^{2}$-distribution.\label{Fig2}}
\end{figure}

For the practical use of a statistic it is important to calculate
or estimate the distribution of the statistic. This can be done by
exact calculation, by approximations and by simulations. Exact calculations
may be both time consuming and difficult. Simulation often requires
statistical insight and programming skills. Therefore most statistical
tests use approximations to calculate the distribution of the statistic.
The distribution of the $\chi^{2}$-statistic becomes closer and closer
to the $\chi^{2}$-distributions as the sample size tends to infinity.
For a large sample size the empirical distribution will with high
probability be close to the generating distribution and any Csisz{\'a}r
$f$-divergence $D_{f}$ can be approximated by a scaled version of
the $\chi^{2}$-statistic
\[
D_{f}\left(P,Q\right)\approx\frac{f^{\prime\prime}\left(0\right)}{2}\cdot\chi^{2}\left(P,Q\right).
\]
Therefore the distribution of any $f$-divergence may be approximated
by a scaled $\chi^{2}$-distribution, i.e. a Gamma distribution. From
this argument one might get the impression that the distribution of
the $\chi^{2}$-statistic is closer to the $\chi^{2}$-distribution.
Figure \ref{Fig1} and Figure \ref{Fig2} show that this is far from
the the case. Figure \ref{Fig1} shows that the distribution of the$G^{2}$-statistic
(i.e. mutual information) is almost as close to a $\chi^{2}$-distribution
as it could be taking into account that mutual information has a discrete
distribution. Each step is intersected very close to its mid point.
Figure \ref{Fig2} shows that the distribution of the Pearson statistic
($\chi^{2}$-statistic) is very different when the observed counts
deviates from the expected distribution and the intersection property
only holds when each cell in the contingency table contains at least
10 observations. These two plots show that at least in some cases
the distribution of mutual information is much closer to a (scaled)
$\chi^{2}$-distribution than the Pearson $\chi^{2}$-statistic is.
The next question is whether there are situations where mutual information
is not approximately $\chi^{2}$-distributed. For contingency tables
that are much less symmetric the intersection property of Figure \ref{Fig1}
is not satisfied when the $G$-statistic is plotted against the $\chi^{2}$-distrbution
so in the rest of this paper a different type of plots will be used.

The use of the $G^{2}$-statistic rather than the $\chi^{2}$-statistic
has become more and more popular since this was recommended in the
1981 edition of the popular textbook of Sokal and Rohlf \cite{Sokal1981}.
A summery of what the typical recommended about whether one should
use the $\chi^{2}$-statistic or the $G^{2}$-Statistic was given
by T. Dunning \cite{Dunning1993}. The short version is that the statistic
is approximately $\chi^{2}$-distributed when each bin contains at
least 5 observations or the calculated variance for each bin is at
least 5, and if any bin contains more than twice the expected number
observations then the $G^{2}$-statistic is preferable to the $\chi^{2}$-statistic.
If the test has only one degree of freedom one often recommend each
bin contain at least 10 observations. 

In this paper we will let $\tau$ denote the circle constant $2\pi$
and let $\phi$ denote the standard Gaussian density
\[
\frac{\exp\left(-\frac{x^{2}}{2}\right)}{\tau^{\nicefrac{1}{2}}}.
\]
 We let $\Phi$ denote the distribution function of the standard Gaussian
\[
\Phi\left(t\right)=\int_{-\infty}^{t}\phi\left(x\right)\,\text{d}x~.
\]

The rest of the paper is organized as follows. In Section \ref{sec:Contingency-tables}
we formulate a general inequality for point probabilities in general
contingency tables in terms of mutual information. In Section \ref{sec:Contingcy2x2}
we intrroduce special notation for $2\times2$ contingency tables,
and we introduce the signed log-likelihood of a $2\times2$ contingency
table. In Section \ref{sec:signedExpo} we introduce the signed log-likelihood
for exponential families. In Section \ref{sec:Inequalities-for-waiting}
we formulate some sharp inequalities for waiting times and as corollaries
we obtain intersection inequalies for binomial distributions and Poisson
distributions. In Section \ref{sec:Numerical-results} we explain
how the intersection property has been checked numerically for all
cases relevant for testing independence. We end with a short discussion.
Most of the proofs, detailed descriptions of the numerical algorithms,
and further plots are given in an appendix. The appendix also contains
some material on larger contingency tables. The appendix can be found
in the arXiv-version of this paper.

\section{Contingency tables\label{sec:Contingency-tables}}

We consider a contingeny table

\begin{center}
\begin{tabular}{|c|c|c|c|c|c|}
\hline 
 &  &  &  &  & Total\tabularnewline
\hline 
 & $X_{11}$ & $X_{12}$ & $\cdots$ & $X_{1\text{\ensuremath{\ell}}}$ & $R_{1}$\tabularnewline
\hline 
 & $X_{21}$ & $X_{22}$ & $\cdots$ & $X_{2\ell}$ & $R_{2}$\tabularnewline
\hline 
 & $\vdots$ & $\vdots$ & $\ddots$ & $\vdots$ & $\vdots$\tabularnewline
\hline 
 & $X_{k1}$ & $X_{k2}$ & $\cdots$ & $X_{k\ell}$ & $R_{k}$\tabularnewline
\hline 
Total & $S_{1}$ & $S_{2}$ & $\cdots$ & $S_{\text{\ensuremath{\ell}}}$ & $N$\tabularnewline
\hline 
\end{tabular} 
\par\end{center}

We fix the row counts $R_{i}$ and the coloum counts $S_{j}$ and
we are interseted in the distribution of the counts $X_{ij}$ under
the null hypothesis that the counted features are independent. These
probabilities are given by
\begin{align*}
\Pr\left(X_{ij}=x_{ij}\right) & =\frac{\left(\begin{array}{c}
N\\
x_{ij}
\end{array}\right)}{\left(\begin{array}{c}
N\\
R_{i}
\end{array}\right)\left(\begin{array}{c}
N\\
S_{j}
\end{array}\right)}\\
 & =\frac{\prod_{i}R_{i}!\prod_{,j}S_{j}!}{N!\prod_{i,j}x_{ij}!}\,.
\end{align*}

The (empirical) mutual information $I$ of the contingency table is
calculated as 
\begin{align*}
N\cdot I & =N\cdot D\left(\left.\hat{P}\right\Vert Q\right)={\textstyle \sum\limits _{i,j}}X_{ij}\ln\frac{x_{ij}}{Nq_{ij}}={\textstyle \sum\limits _{i,j}}X_{ij}\ln\frac{N\cdot x_{ij}}{R_{i}\cdot S_{j}}\\
 & =\sum_{i,j}x_{ij}\ln\left(x_{ij}\right)-\sum_{i}R_{i}\ln\left(R_{i}\right)-\sum_{j}S_{j}\ln\left(S_{j}\right)\\
 & \,\,+N\ln\left(N\right).
\end{align*}

\begin{thm}
\label{thm:For-any-contingency}For a contingency table without empty
cells the following inequality holds
\begin{equation}
\Pr\left(X_{ij}=x_{ij}\right)\leq\mathrm{e}^{^{-N\cdot I}}\left(\frac{\prod_{i}R_{i}\prod_{,j}S_{j}}{\tau^{\left(k-1\right)\left(\ell-1\right)}N\prod_{i,j}x_{ij}}\right)^{\nicefrac{1}{2}}.\label{eq:GenerelIneq}
\end{equation}

\end{thm}
Although Inequality \ref{eq:GenerelIneq} is quite tigth for most
values of $x$ it cannot be used to provide bounds on the tail probabilities
because it is not valid when one of the cells is empty. Inequalities
\[
\Pr\left(X_{ij}=x_{ij}\right)\geq C\mathrm{e}^{^{-N\cdot I}}\left(\frac{\prod_{i}R_{i}\prod_{,j}S_{j}}{\tau^{\left(k-1\right)\left(\ell-1\right)}N\prod_{i,j}x_{ij}}\right)^{\nicefrac{1}{2}}
\]
where $C<1$ is a positive constant can also be proved, but the optimal
constant will depend on the sample size and on the size of the contingency
table.

\section{Contingency tables with one degree of freedom\label{sec:Contingcy2x2}}

In this paper our theoretical results will focus on contingency tables
with one degree of freedom and there are several reasons for this.
The distribution of the statistic is closely related to the hypergeometric
distribution that is well studied in the literature. It is know that
the use of the $\chi^{2}$-distributions is most problematic when
the number of degrees of freedom is small. Results of testing independency
can be compared with similar results for testing goodness of fit that
are only available for one degree of freedom.

First we will introduce some special notation. We consider a contingency
table

\begin{center}
\begin{tabular}{|c|c|c|c|}
\hline 
 &  &  & Total\tabularnewline
\hline 
 & $X$ & $r-X$ & $r$\tabularnewline
\hline 
 & $n-X$ & $N-r+X$ & $N-r$\tabularnewline
\hline 
Total & $n$ & $N-n$ & $N$\tabularnewline
\hline 
\end{tabular}
\par\end{center}

\begin{figure}
\centering{}\definecolor{ffqqqq}{rgb}{1,0,0} 
\definecolor{qqqqff}{rgb}{0,0,1} \begin{tikzpicture}[line cap=round,line join=round,>=triangle 45,x=0.35714285714285715cm,y=0.35714285714285715cm] 
\draw[->,color=black] (-5.5,0) -- (8.5,0); 
\foreach \x in {-4,-2,2,4,6,8} 
\draw[shift={(\x,0)},color=black] (0pt,2pt) -- (0pt,-2pt) node[below] {\footnotesize $\x$}; 
\draw[->,color=black] (0,-5.5) -- (0,8.5); 
\foreach \y in {-4,-2,2,4,6,8} 
\draw[shift={(0,\y)},color=black] (2pt,0pt) -- (-2pt,0pt) node[left] {\footnotesize $\y$}; 
\draw[color=black] (0pt,-10pt) node[right] {\footnotesize $0$}; 
\clip(-5.5,-5.5) rectangle (8.5,8.5); 
\draw [color=qqqqff] (-3.3586706055,-2.9212013103)-- (-3.3586706055,-3.7711109175); 
\draw [color=qqqqff] (-2.5548926947,-2.161775638)-- (-2.5548926947,-2.9212013103); 
\draw [color=qqqqff] (-1.8153892123,-1.445635242)-- (-1.8153892123,-2.161775638); 
\draw [color=qqqqff] (-1.1095557208,-0.7540939007)-- (-1.1095557208,-1.445635242); 
\draw [color=qqqqff] (-0.4231437738,-0.0767183138)-- (-0.4231437738,-0.7540939007); 
\draw [color=qqqqff] (0.924060006,1.2626508149)-- (0.924060006,0.5936264999); 
\draw [color=qqqqff] (1.596776096,1.9356064705)-- (1.596776096,1.2626508149); 
\draw [color=qqqqff] (2.2761178875,2.6179906906)-- (2.2761178875,1.9356064705); 
\draw [color=qqqqff] (2.9680823484,3.3163109937)-- (2.9680823484,2.6179906906); 
\draw [color=qqqqff] (3.6802323215,4.0392746932)-- (3.6802323215,3.3163109937); 
\draw [color=qqqqff] (4.4234304316,4.8002254118)-- (4.4234304316,4.0392746932); 
\draw [color=qqqqff] (5.2159712851,5.6236661561)-- (5.2159712851,4.8002254118); 
\draw [color=qqqqff] (6.0975364584,6.5717816881)-- (6.0975364584,5.6236661561); 
\draw [color=qqqqff] (0.2525852384,0.5936264999)-- (0.2525852384,-0.0767183138); 
\draw [color=qqqqff] (7.27,6.57)-- (7.27,8.5); 
\draw [color=qqqqff] (0.2525852384,0.5936264999)-- (0.2525852384,-0.0767183138); 
\draw [color=qqqqff] (-4.39,-3.77)-- (-4.39,-5.5); 
\draw[color=ffqqqq,smooth,samples=100,domain=-5.5:8.5] plot(\x,{(\x)}); \end{tikzpicture}\caption{\label{fig:QuantilesHypVsGaus201010}Quantiles of $g\left(X\right)$
against a standard Gaussian when $X$ is hypergeometric with parameters
$N=40,$ $n=15$, and $r=15.$ The streight line corresponds to a
perfect match with a Gaussian distribution.}
\end{figure}
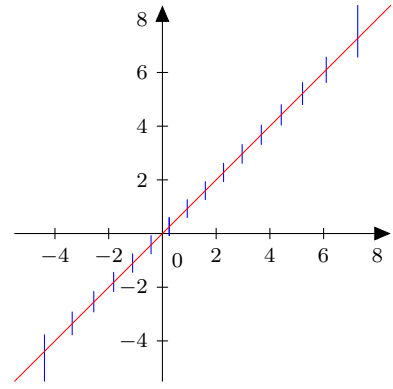

If the marginal counts are kept fixed and independence is assumed
then the distribution of $X$ is a hypergeometric distribution with
parameters $\left(N,r,n\right)$ and point probabilities 
\[
P\left(X=x\right)=\frac{\left(\begin{array}{c}
r\\
x
\end{array}\right)\left(\begin{array}{c}
N-r\\
n-x
\end{array}\right)}{\left(\begin{array}{c}
N\\
n
\end{array}\right)}
\]
This hypergeometric disribution has mean value $\frac{nr}{N}$ and
variance $\frac{nr\left(N-r\right)\left(N-n\right)}{N^{2}\left(N-1\right)}$.
For $2\times2$ contingency tables we introduce a random variable
$G$ that is related to mutual information in the same way as $\chi$
is related to $\chi^{2}$.
\begin{defn}
The \emph{signed log-likelihood} of a $2\times2$ contingency table
is defined as
\[
G\left(X\right)=\left\{ \begin{array}{cc}
-\left(2N\cdot I\left(X\right)\right)^{\nicefrac{1}{2}}, & \text{for }X<\frac{nr}{N};\\
+\left(2N\cdot I\left(X\right)\right)^{\nicefrac{1}{2}}, & \text{for }X\geq\frac{nr}{N}.
\end{array}\right.
\]

\end{defn}
With this notation the first factor in Inequality \ref{eq:GenerelIneq}
can be written as $\frac{\mathrm{e}^{^{-N\cdot I}}}{\tau^{\nicefrac{1}{2}}}=\phi\left(G\left(x\right)\right).$
A QQ-plot of $G\left(X\right)$ against a standard Gaussian is given
in Figure \ref{fig:QuantilesHypVsGaus201010}. This plot and other
similar plots support the following \emph{intersection conjecture
for hypergeometric distributions}.
\begin{conjecture}
\label{hyperIntersectionConjecture}For a $2\times2$ contingency
table let $G\left(X\right)$ denote the signed log-likelihod. Then
the quantile transform between $G\left(X\right)$ and a standard Gaussian
$Z$ is a step function and the identity function intersects each
step, i.e. $\Pr\left(X<x\right)<\Pr\left(Z\leq g\left(x\right)\right)<\Pr\left(X\leq x\right)$
for all integers $x.$ 
\end{conjecture}
If $n$ is much smaller than $N$ then the distribution of $X$ can
be approximated by binomial distribution $bin\left(n,p\right)$ where
$p=\nicefrac{r}{N}$. If $n$ is much smaller than $N$ we can approximate
the mutual information by the information divergence $\frac{x}{n}\ln\frac{\nicefrac{x}{n}}{p}+\left(1-\frac{x}{n}\right)\ln\left(\frac{1-\nicefrac{x}{n}}{1-p}\right).$
Therefore it is of interest to see if the intersection conjecture
holds in the limiting case of a binomial distribution. 

\section{The signed log-likelihood for exponential families\label{sec:signedExpo}}

Consider the 1-dimensional exponential family $P_{\beta}$ where
\[
\frac{\text{d}P_{\beta}}{\text{d}P_{0}}\left(x\right)=\frac{\exp\left(\beta\cdot x\right)}{Z\left(\beta\right)}
\]
 and $Z$ denotes the \emph{moment generating function} given by Z$\left(\beta\right)=\int\exp\left(\beta\cdot x\right)\,\text{d}P_{0}x.$
Let $P^{\mu}$ denote the element in the exponential family with mean
value $\mu,$ and let $\hat{\beta}\left(\mu\right)$ denote the corresponding
maximum likelihood estimate of $\beta.$ Let $\mu_{0}$ denote the
mean value of $P_{0}.$ Then
\begin{align*}
D\left(P^{\mu}\Vert P_{0}\right) & =\int\ln\left(\frac{\text{d}P^{\mu}}{\text{d}P_{0}}\left(x\right)\right)\,\text{d}P^{\mu}x.\\
 & =\mu\cdot\hat{\beta}\left(\mu\right)-\ln\left(Z\left(\hat{\beta}\left(\mu\right)\right)\right)
\end{align*}

\begin{defn}
Let $X$ be a random variable with distribution $P_{0}.$ Then the
$G$-transform $g\left(X\right)$ of $X$ is the random variable given
by
\[
G\left(x\right)=\left\{ \begin{array}{cc}
-\left(2D\left(P^{x}\Vert P_{0}\right)\right)^{\nicefrac{1}{2}}, & \text{for }x<\mu_{0};\\
+\left(2D\left(P^{x}\Vert P_{0}\right)\right)^{\nicefrac{1}{2}}, & \text{for }x\geq\mu_{0}.
\end{array}\right.
\]
\end{defn}
\begin{thm}
Let $X$ denote a binomial distributed or Poisson distributed random
variable and let $G\left(X\right)$ denote the $G$-transform of $X.$
The quantile transform between $g\left(X\right)$ and a standard Gaussian
$Z$ is a step function and the identity function intersects each
step, i.e. $P\left(X<k\right)<P\left(Z\leq g\left(k\right)\right)<P\left(X\leq k\right)$
for all integers $k.$ 
\end{thm}
For Poisson distributions this result was proved in \cite{Harremoes2012}
where the result for binomial distributions was also conjectured.
A proof of this conjecture was recently given by Serov and Zubkov
\cite{Zubkov2013}. Were we will prove that the intersection results
for binomial distributions and Poisson distributions follows from
more general results for waiting times. We will need the following
general lemma.
\begin{lem}
\label{lem:exponentialDiff}In an exponential family with $G=\pm\left(2D\left(\left.P^{\mu}\right\Vert P^{\mu_{0}}\right)\right)^{\nicefrac{1}{2}}$
one has 
\[
\frac{\partial}{\partial\mu_{0}}\Phi\left(G\right)=\frac{\phi\left(G\right)}{V\left(\mu_{0}\right)}\cdot\frac{\mu_{0}-\mu}{G}.
\]

\end{lem}

\section{Inequalities for waiting times\label{sec:Inequalities-for-waiting}}

Next we will formulate some inequalities for waiting times. Let $nb\left(p,\ell\right)$
denote distribution of the number of failures before the $\ell$'th
succes in a Bernoulli process with success probability $p.$ The mean
value is given by $\mu=\ell\frac{1-p}{p}$. 
\begin{lem}
\label{lem:difNeg}If the distribution if $W_{1}$ is $nb\left(p,\ell\right)$
with $\mu=\ell\frac{1-p}{p}$ then the partial derivative of the point
probaiblity equals
\[
\frac{\partial}{\partial\mu}\Pr\left(W_{1}=k\right)=\Pr\left(W_{2}=k-1\right)-\Pr\left(W_{2}=k\right).
\]
where $W_{2}$ is $nb\left(p,\ell+1\right).$\end{lem}
\begin{thm}
\label{thm:Negbin}For any $\ell\in\left[1;\infty\right[$ the random
variable $W_{2}$ with distribution $Nb\left(p,\ell\right)$satisfies
\[
\Phi\left(G_{nb\left(p,\ell\right)}\left(k\right)\right)\leq\Pr\left(W_{2}\leq k\right),
\]
where $G_{nb\left(p,\ell\right)denotes}$ the signed log-likelihood
of the negative binomial distribution.\end{thm}
\begin{IEEEproof}
For any $k,l$ one has to prove that $\Pr\left(W_{2}\leq k\right)-\Phi\left(G_{nb\left(p,\ell\right)}\left(k\right)\right)$
is non-negative. If $p=0$ or $p=1$ this is obvious. The result is
obtained by differentiating several times with respect to the mean
value. See Appendix for details. 
\end{IEEEproof}
A Gamma distribution can be considered as a limit of negative binomial
distributions, which leads to the following corollary.
\begin{cor}
If $W$ denotes a Gamma distribution then
\[
\Phi\left(G_{\Gamma}\left(x\right)\right)\leq\Pr\left(W\leq x\right).
\]

\end{cor}
A similar inequality is satisfied by inverse Gaussian distributions
that are waiting times for browning motions.
\begin{thm}
\label{thm:InverseGaussian}For any $x\in\left[0;\infty\right[$ the
random variable $W$ with distribution $IG\left(\mu,\lambda\right)$
satisfies
\[
\Phi\left(G_{IG\left(\mu,\lambda\right)}\left(k\right)\right)\leq\Pr\left(W\leq k\right),
\]
where $G_{IG\left(\mu,\lambda\right)}$ denotes the signed log-likelihood
of the inverse Gaussian.\end{thm}
\begin{IEEEproof}
The proof uses the same technique as in \cite{Harremoes2012}. See
Appendix.
\end{IEEEproof}
Assume that $M$ is binomial $bin\left(n,p\right)$ and $W$ is negative
binomial $nb\left(p,\ell\right).$ Then
\[
\Pr\left(M\geq\ell\right)=\Pr\left(W+\ell\leq n\right).
\]
The divergence can be calcuated as
\begin{multline*}
D\left(\left.nb\left(p_{1},\ell\right)\right\Vert nb\left(p_{2}\ell\right)\right)=\\
\frac{\ell}{p_{1}}\left(p_{1}\ln\frac{p_{1}}{p_{2}}+\left(1-p_{1}\right)\ln\frac{1-p_{1}}{1-p_{2}}\right).
\end{multline*}
 If $\frac{\ell}{p_{1}}=n$ then
\begin{multline*}
D\left(\left.nb\left(p_{1},\ell\right)\right\Vert nb\left(p_{2}\ell\right)\right)\\
=n\left(p_{1}\ln\frac{p_{1}}{p_{2}}+\left(1-p_{1}\right)\ln\frac{1-p_{1}}{1-p_{2}}\right)\\
=D\left(\left.bin\left(n,\frac{\ell}{n}\right)\right\Vert bin\left(n,p\right)\right).
\end{multline*}
If $G_{bin}$ is the signed likelihood log-likelihood of $bin\left(n,p\right)$
and $G_{nb}$is the signed log-likelihood of $nb\left(p,\ell\right)$
then $G_{bin}\left(\ell\right)=-G_{nb}\left(n\right).$ This shows
that if the log-likelihood of a the negative binomial distributions
are close to Gaussian then the same is true for the signed log-likelihood
of binomial distribtions. Our inequality for the negative binomial
distribution can be translated into an inequality for the binomial
distribution.
\begin{thm}
\label{thm:Binomialintersection}Binomial distributions satisfy the
intersection property.\end{thm}
\begin{IEEEproof}
We have to prove that
\[
\Pr\left(S_{n}<\ell\right)\leq\Phi\left(G_{bin\left(n,p\right)}\left(\ell\right)\right)\leq\Pr\left(S_{n}\leq\ell\right).
\]
We will only prove the left inequality since the right inequality
follows from the left inequality by replacing $p$ by $1-p$ and replacing
$\ell$ by $n-\ell.$ We have
\begin{eqnarray*}
\Pr\left(S_{n}<\ell\right) & = & 1-\Pr\left(S_{n}\geq\ell\right)=1-\Pr\left(W_{1}\leq n\right)\\
 & \leq & 1-\Phi\left(G_{nb\left(p,\ell\right)}\left(n\right)\right)=\Phi\left(-G_{nb\left(p,\ell\right)}\left(n\right)\right)\\
 & = & \Phi\left(G_{bin\left(n,p\right)}\left(\ell\right)\right).
\end{eqnarray*}

Note that Theorem \ref{thm:Negbin} cannot be proved from Theorem
\ref{thm:Binomialintersection} because the number parameter for a
binomial distribution has to be an integer while the number parameter
of a negative binomial distribution may assume any positive value.Since
a Poisson distribution is a limit of binomial distributions we get
the following corollary.\end{IEEEproof}
\begin{cor}
Poisson distributions have the intersection property.
\end{cor}
These results allow us to sharpen a bound for the tail probabilities
of a hypergeometric distribution proved by Chv{\'a}tal \cite{Chvatal1979}.
\begin{cor}
Let the distribution of $X$ be hypergeometric with parameters $\left(N,n,r\right)$
and assume that $x\leq nr/N.$ Then 
\begin{equation}
\Pr\left(X<\ell\right)<\Phi\left(-\left(2nD\left(\left.\left(\frac{x}{n},1-\frac{x}{n}\right)\right\Vert \left(\frac{r}{N},1-\frac{r}{N}\right)\right)\right)\right).\label{eq:ChvaImproved}
\end{equation}

\end{cor}
If we assume that mutual information has (scaled) $\chi^{^{2}}$-distribution
then we get an estimate of $\Pr\left(X<\ell\right)$ that is less
than $\Phi\left(-\left(2nD\left(\left.\left(\frac{x}{n},1-\frac{x}{n}\right)\right\Vert \left(\frac{r}{N},1-\frac{r}{N}\right)\right)\right)\right)$
but if we assume that the $\chi^{2}$-statistic is $\chi^{2}$-distributed
we often get an estimate of $\Pr\left(X<\ell\right)$ that violates
Inequality \ref{eq:ChvaImproved} and in this sense mutual information
has a distribution that is closer to the $\chi^{2}$-distribution
than the distribution of the $\chi^{2}$-statistic is.

\section{Numerical results\label{sec:Numerical-results}}

We have not yet been able to prove the intersection conjecture \ref{hyperIntersectionConjecture}
for all hypergeometric distributions, but we have numerical calculations
that confirm the conjecture for all cases that are relevant for testing
independence. 

Due to a symmetry of the problem we only have to check that $\Pr\left(X<x\right)<\Pr\left(Z\leq g\left(x\right)\right)$
holds for all $2\times2$ contingency tables. For any finite sample
size there are only finitely many contingency table with this sample
size, so for any fixed sample size the inequality can be checked numerically.
We have checked the hypergeometric intersection conjecture for any
sample size up to 200. 

If all bins contain at least 10 observations then the rule of thumb
states that the $\chi^{2}$-statistic and the $G^{2}$-statistic will
both follow a $\chi^{2}$-distribution, so it should be sufficient
to test the hypothesis for $x<10.$

As a rule of thumb the approximation by a binomial distribution is
applicable when $n\leq0.1\cdot N$ and the intersection property holds
for binomial distributions. Therefore we are only interested in the
cases where $n>0.1\cdot N$. By symmetry in the parameters it is also
possible to approximate the distribution of $X$ by the binomial distribution
$bin\left(r,\nicefrac{n}{N}\right)$ if $r\leq0.1\cdot N$.

For testing independence one will normally choose a significance level
$\alpha\in\left[0.001,0.1\right]$. If the $p$-value can easily be
bounded away from this interval there is no need to make a refined
estimate of the $p$-value via the intersection conjecture for the
hypergeometric distribution. For $x\leq9$ the probability $\Pr\left(X\leq x\right)$
is bounded by $\Pr\left(X\leq9\right).$ Let $\lambda$ denote the
mean of $X$. Then $\Pr\left(X\leq9\right)\leq\Pr\left(T\leq9\right)$
where $T$ is Poisson distributed with mean $\lambda$. Now $\Pr\left(T\leq9\right)\leq0.000974$
which is less that a significance level of 0.001 for $\lambda\geq22.7$
so we are most interested in the cases where $\lambda\leq22.7$. If
we combine this with the condition $n\geq0.1\cdot N$ and $r\geq0.1\cdot N$
and that the mean of $X$ is $\frac{nr}{N}$ we see that we only have
to check cases where $N\leq100\cdot22.7=2270$. If we impose all these
conditions we are left with a high but finite number of cases that
can be checked numerically one by one. Further details are given in
the appendix. The result is that the intersection conjecture holds
for hypergeometric distributions for all cases that are of interest
in relation to testing independence.

\section{Discussion\label{sec:Discussion}}

The distribution of the signed log-likelihood is close a standard
Gaussian for a variety of distributions. As an asymptotic result for
large sample sizes this is not new \cite{Zhang2008}, but for the
most important distributions like the binomial distributions, the
Poisson distributions, the negative binomial distributions, the inverse
Gaussian distributions and the Gamma distributions we can formulate
sharp inequalities that hold for any sample size. All these distributions
have variance functions that are polynomials of order 2 and 3. Natural
exponential families with polynomial variance functions of order at
most 3 have been classified \cite{Morris1982,Letac1990} and there
is a chance that one can formulate and prove a sharp inequality for
each of these exponential families.

As we have seen there is strong numerical evidence that an intersection
conjecture holds for all hypergeometric distributions, but a proof
will be more difficult because the hypergeometric distributions do
not form an exponential family and no variance function is available.
In order to prove the conjecture it may prove useful to note that
any hypergeometric distribution is a Poisson binomial distribution
\cite{Barbour1992}, i.e. the distribution of a sum of independent
Bernoulli random variables with different success probabilities. Exploration
of intersection conjectures for Poisson binomial distributions is
work in progress and seems to be the most promising method direction
in order to verify the intersection conjecture for the hypergeometric
distributions.

\section{Acknowledgment}

The author want to thank Gabor Tusn{\'a}dy, L{\'a}szlo Gy{\"o}rfi, Joram Gat,
Jan{\'o}s Komlos, and A. Zubkov for useful discussions.

\bibliographystyle{ieeetr}

\onecolumn

\appendix
This appendix contain proofs and further material that supports the
conclusion of the paper.

\section*{Proof of Theorem \ref{thm:For-any-contingency}}

For the proof of our next theorem we need the following lemma. 
\begin{lem}
\label{lem:Lagrange}The function $\prod_{i=1}^{n}\left(1+\frac{1}{x_{i}}\right)$
is minimal under the constraint $\sum_{i=1}^{n}x_{i}=\lambda$ and
$x_{i}>0$ all $x_{i}$ are equal.\end{lem}
\begin{IEEEproof}
We use the method of Lagrange multipliers. The matrix of partial derivatives
is 
\[
\left(\begin{array}{cccc}
-\frac{1}{x_{1}^{2}}\cdot\frac{\prod_{i=1}^{n}\left(1+\frac{1}{x_{i}}\right)}{1+\frac{1}{x_{1}}} & -\frac{1}{x_{2}^{2}}\cdot\frac{\prod_{i=1}^{n}\left(1+\frac{1}{x_{i}}\right)}{1+\frac{1}{x_{2}}} & \cdots & -\frac{1}{x_{n}^{2}}\cdot\frac{\prod_{i=1}^{n}\left(1+\frac{1}{x_{i}}\right)}{1+\frac{1}{x_{n}}}\\
1 & 1 & \cdots & 1
\end{array}\right).
\]
This matrix is singular if $x_{i}^{2}\left(1+\frac{1}{x_{i}}\right)$
is constant, but $x_{i}^{2}\left(1+\frac{1}{x_{i}}\right)=x^{2}+x$
which is an increasing function. Hence all $x_{i}$ must be equal.
\end{IEEEproof}
The upper bound requires a detailed calculation where we use lower
and upper bounds related to Stirling's approximations of factorials.
\[
n^{n}\mathrm{e}^{-n}\left(\tau n\right)^{1/2}\left(1+\frac{1}{12n+1}\right)\leq n!\leq n^{n}\mathrm{e}^{-n}\left(\tau n\right)^{1/2}\left(1+\frac{1}{12n}\right).
\]

We have 
\begin{align*}
\Pr\left(X_{ij}=x_{ij}\right) & =\frac{\prod_{i}R_{i}!\prod_{,j}S_{j}!}{N!\prod_{i,j}x_{ij}!}\\
 & \leq\frac{\prod_{i}\left(R_{i}^{R_{i}}\exp\left(-R_{i}\right)\left(\tau R_{i}\right)^{\nicefrac{1}{2}}\left(1+\frac{1}{12R_{i}}\right)\right)\prod_{,j}\left(S_{j}^{S_{j}}\exp\left(-S_{j}\right)\left(\tau S_{j}\right)^{\nicefrac{1}{2}}\left(1+\frac{1}{12S_{j}}\right)\right)}{\left(N^{N}\exp\left(-N\right)\left(\tau N\right)^{\nicefrac{1}{2}}\left(1+\frac{1}{12N+1}\right)\right)\prod_{i,j}\left(X_{ij}^{X_{ij}}\exp\left(-x_{ij}\right)\left(\tau x_{ij}\right)^{\nicefrac{1}{2}}\left(1+\frac{1}{12x_{ij}+1}\right)\right)}\\
 & =\frac{\prod_{i}R_{i}^{R_{i}}\prod_{,j}\left(S\right)}{N\prod_{i,j}X_{ij}^{X_{ij}}}\cdot\frac{\tau^{\nicefrac{k}{2}}\cdot\tau^{\nicefrac{\ell}{2}}}{\tau^{\nicefrac{1}{2}}\cdot\tau^{\nicefrac{k\ell}{2}}}\left(\frac{\prod_{i}R_{i}\prod_{,j}S_{j}}{N\prod_{i,j}x_{ij}}\right)^{\nicefrac{1}{2}}\cdot\frac{\prod_{i}\left(1+\frac{1}{12R_{i}}\right)\prod_{,j}\left(1+\frac{1}{12S_{j}}\right)}{\left(1+\frac{1}{12N+1}\right)\prod_{i,j}\left(1+\frac{1}{12x_{ij}+1}\right)}
\end{align*}
The first factor equals $\exp\left(-N\cdot I\right)$ so it is sufficient
to prove that the last factor is less than 1. We have 
\[
\prod_{i,j}\left(1+\frac{1}{12x_{ij}+1}\right)=\prod_{i}\left(\prod_{j}\left(1+\frac{1}{12x_{ij}+1}\right)\right)^{\nicefrac{1}{2}}\cdot\prod_{j}\left(\prod_{i}\left(1+\frac{1}{12x_{ij}+1}\right)\right)^{\nicefrac{1}{2}}
\]
Now 
\begin{align*}
\left(\prod_{j}\left(1+\frac{1}{12x_{ij}+1}\right)\right)^{\nicefrac{1}{2}} & \geq\left(\prod_{j}\left(1+\frac{1}{\frac{12R_{i}}{\text{\ensuremath{\ell}}}+1}\right)\right)^{\nicefrac{1}{2}}\\
 & =\left(1+\frac{1}{\frac{12R_{i}}{\text{\ensuremath{\ell}}}+1}\right)^{\nicefrac{\ell}{2}}.
\end{align*}
We now have to prove that 
\[
\left(1+\frac{1}{\frac{12R_{i}}{\text{\ensuremath{\ell}}}+1}\right)^{\nicefrac{\ell}{2}}\geq1+\frac{1}{12R_{i}}
\]
for $\text{\ensuremath{\ell\geq}2.}$ For $\ell=2$ the inequality
is obviously satisfied because we have assumed that none of the cells
were empty. It is sufficient to prove that 
\[
\left(1+\frac{1}{\frac{12R_{i}}{\text{\ensuremath{\ell}}}+1}\right)^{\nicefrac{\ell}{2}}=\exp\left(\frac{\ell}{2}\ln\frac{12R+2\ell}{12R+\ell}\right)
\]
We see that $\frac{\ell}{2}\ln\frac{12R+2\ell}{12R+\ell}$ is a product
of two positive increasing functions and therefore it is increasing.

\section*{The $\chi^{2}$-statistic for $2\times2$ contingency tables}

\begin{table}
\centering{}%
\begin{tabular}{|c|c|c|c|}
\hline 
 &  &  & Total\tabularnewline
\hline 
\hline 
 & $X-\frac{nr}{N}$ & $\frac{nr}{N}-X$ & $0$\tabularnewline
\hline 
 & $\frac{nr}{N}-X$ & $X-\frac{nr}{N}$ & $0$\tabularnewline
\hline 
Total & $0$ & $0$ & $0$\tabularnewline
\hline 
\end{tabular}\caption{\label{tab:Deviations}Deviations from the expected values.}
\end{table}
Then the table of differences between observed and expected values
can be seen in Table \ref{tab:Deviations} and the 
\begin{align*}
\chi^{2} & =\frac{\left(X-\frac{nr}{N}\right)^{2}}{\frac{nr}{N}}+\frac{\left(\frac{nr}{N}-X\right)^{2}}{\frac{n\left(N-r\right)}{N}}+\frac{\left(\frac{nr}{N}-X\right)^{2}}{\frac{r\left(N-n\right)}{N}}+\frac{\left(X-\frac{nr}{N}\right)^{2}}{\frac{\left(N-r\right)\left(N-r\right)}{N}}\\
 & =\frac{\left(X-\frac{nr}{N}\right)^{2}}{\frac{nr\left(N-r\right)\left(N-n\right)}{N^{3}}}
\end{align*}
If we define $\chi=\pm\left(\chi^{2}\right)^{1/2}$ then we see that
\[
\chi=\frac{X-\frac{nr}{N}}{\left(\frac{nr\left(N-r\right)\left(N-n\right)}{N^{3}}\right)^{\nicefrac{1}{2}}}
\]

This hypergeometric distribution has mean value $\frac{nr}{N}$ and
variance $\frac{nr\left(N-r\right)\left(N-n\right)}{N^{2}\left(N-1\right)}$.
Hence 
\[
\chi=\frac{X-E\left[X\right]}{\sigma\left(X\right)}\cdot\left(\frac{N}{N-1}\right)^{\nicefrac{1}{2}}
\]
Since $X$ is approximately Gaussian the distribution of the first
factor tends to a standard Gaussian and the second factor tends to
1, so the $\chi^{2}$-statistic is approximately $\chi^{2}$-distributed.

\section*{Proof of Lemma \ref{lem:exponentialDiff}}

From $G^{2}=2D$ follows that $G\cdot G'=D'$ and that
\[
\frac{\partial G}{\partial\mu_{0}}=\frac{\frac{\partial D}{\partial\mu_{0}}}{G}.
\]
Hence
\begin{eqnarray*}
\frac{\partial}{\partial\mu_{0}}\Phi\left(G\right) & = & \phi\left(G\right)\cdot\frac{\partial G}{\partial\mu_{0}}\\
 & = & \frac{\phi\left(G\right)}{G}\cdot\frac{\partial D}{\partial\mu_{0}}.
\end{eqnarray*}
We write the divergence as $D=\beta_{\mu_{0}}\cdot\mu-\ln Z\left(\beta_{\mu_{0}}\right).$
Hence
\begin{eqnarray*}
\frac{\partial D}{\partial\mu_{0}} & = & -\frac{\partial\beta_{\mu_{0}}}{\partial\mu_{0}}\cdot\mu+\frac{Z'\left(\beta_{\mu_{0}}\right)}{Z\left(\beta_{\mu_{0}}\right)}\cdot\frac{\partial\beta_{\mu_{0}}}{\partial\mu_{0}}\\
 & = & \left(\frac{Z'\left(\beta_{\mu_{0}}\right)}{Z\left(\beta_{\mu_{0}}\right)}-\mu\right)\frac{1}{\frac{\partial\mu_{0}}{\partial\beta_{\mu_{0}}}}\\
 & = & \frac{\mu_{0}-\mu}{V\left(\mu_{0}\right)}
\end{eqnarray*}
and the result follows.

\section*{Proof of Lemma \ref{lem:difNeg}}

We have 
\begin{eqnarray*}
\frac{\partial}{\partial\mu}\Pr\left(W_{1}=k\right) & = & \frac{1}{\frac{d\mu}{dp}}\cdot\frac{\partial}{\partial p}\left(\frac{\ell^{\bar{k}}}{k!}p^{\ell}\left(1-p\right)^{k}\right)\\
 & = & \frac{1}{-\ell p^{-2}}\left(\frac{\ell^{\bar{k}}}{k!}\ell p^{\ell-1}\left(1-p\right)^{n}-\frac{\ell^{\bar{k}}}{k!}p^{\ell}k\left(1-p\right)^{k-1}\right)\\
 & = & -\frac{\ell^{\bar{k}}}{k!}p^{\ell+1}\left(1-p\right)^{k}+\frac{\ell^{\bar{k}}}{\ell\left(k-1\right)!}p^{\ell+2}\left(1-p\right)^{k-1}\\
 & = & \frac{\left(\ell+1\right)^{\bar{k-1}}}{\left(k-1\right)!}p^{\ell+2}\left(1-p\right)^{k-1}-\frac{\left(\ell+1\right)^{\bar{k}}-k\cdot\left(\ell+1\right)^{\bar{k-1}}}{k!}p^{\ell+1}\left(1-p\right)^{k}\\
 & = & \left(\frac{\left(\ell+1\right)^{\bar{k-1}}}{\left(k-1\right)!}p^{\ell+2}\left(1-p\right)^{k-1}+\frac{\left(\ell+1\right)^{\bar{k-1}}}{\left(k-1\right)!}p^{\ell+1}\left(1-p\right)^{k}\right)-\frac{\left(\ell+1\right)^{\bar{k}}}{k!}p^{\ell+1}\left(1-p\right)^{k}\\
 & = & \frac{\left(\ell+1\right)^{\bar{k}-1}}{\left(k-1\right)!}p^{\ell+1}\left(1-p\right)^{k-1}-\frac{\left(\ell+1\right)^{\bar{k}}}{k!}p^{\ell+1}\left(1-p\right)^{k}\\
 & = & \Pr\left(W_{2}=k-1\right)-\Pr\left(W_{2}=k\right).
\end{eqnarray*}

As a consequence 
\begin{eqnarray*}
\frac{\partial}{\partial\mu}\Pr\left(W_{1}\leq k\right) & = & -\Pr\left(W_{2}=k\right)\\
 & = & -\frac{\left(\ell+1\right)^{\bar{k}}}{k!}p^{\ell+1}\left(1-p\right)^{k}.
\end{eqnarray*}

\section*{Proof of Theorem \ref{thm:Negbin}}

Introduce
\[
\delta\left(\mu_{0}\right)=\Pr\left(W\leq k\right)-\Phi\left(G\left(k\right)\right)
\]
and note that $\delta\left(0\right)=\lim_{\mu_{0}\to\infty}\delta\left(\mu\right)=0.$
It is sufficient to prove that $\delta$ is first increasing and then
decreasing in $\left[0,\infty\right[.$

According to Lemma \ref{lem:difNeg} the derivative of $\Pr\left(W\leq k\right)$
with respect to $\mu_{0}$ is 
\begin{align*}
\frac{\partial}{\partial\mu_{0}}\Pr\left(W\leq k\right) & =-\frac{\left(\ell+1\right)^{\bar{k}}}{k!}p_{0}^{\ell+1}\left(1-p_{0}\right)^{k}.
\end{align*}
Put $\hat{p}=\frac{\ell}{k+\ell}$ and
\begin{align*}
D & =D\left(\left.nb\left(\hat{p},\ell\right)\right\Vert nb\left(p_{0},\ell\right)\right)\\
 & =nD\left(\left.\left(\hat{p},1-\hat{p}\right)\right\Vert \left(p_{0},1-p_{0}\right)\right)
\end{align*}
and $G=sign\left(p_{0}-\hat{p}\right)\cdot\left(2D\right)^{\nicefrac{1}{2}}$.
We use that
\begin{align*}
p_{0}^{\ell}\left(1-p_{0}\right)^{k} & =\exp\left(\ell\ln\left(p_{0}\right)+k\ln\left(1-p_{0}\right)\right)\\
 & =\exp\left(-D-\left(k+\ell\right)H\left(\hat{p},1-\hat{p}\right)\right)
\end{align*}
so that
\begin{align*}
\frac{\partial}{\partial\mu_{0}}\Pr\left(W\leq k\right) & =-\frac{\left(\ell+1\right)^{\bar{k}}}{k!}p_{0}^{\ell+1}\left(1-p_{0}\right)^{k}\\
 & =-\frac{\left(\ell+1\right)^{\bar{k}}p_{0}}{k!}\exp\left(-D-\left(k+\ell\right)H\left(\hat{p},1-\hat{p}\right)\right)\\
 & =-\frac{\left(\ell+1\right)^{\bar{k}}\tau^{\nicefrac{1}{2}}p_{0}}{k!\exp\left(nH\left(\hat{p},1-\hat{p}\right)\right)}\cdot\phi\left(G\right).
\end{align*}
Lemma \ref{lem:exponentialDiff} implies that
\begin{align*}
\frac{\partial}{\partial\mu_{0}}\Phi\left(G\right) & =\frac{\phi\left(G\right)}{V}\cdot\frac{\mu_{0}-\mu}{G}\\
 & =\phi\left(G\right)p_{0}\cdot\frac{\mu_{0}-\mu}{\mu_{0}\cdot G}.
\end{align*}
Combining these results we get
\begin{align*}
\frac{\partial\delta}{\partial\mu_{0}} & =-\frac{\left(\ell+1\right)^{\bar{k}}\tau^{\nicefrac{1}{2}}p_{0}}{k!\exp\left(nH\left(\hat{p},1-\hat{p}\right)\right)}\cdot\phi\left(G\right)-\phi\left(G\right)p_{0}\cdot\frac{\mu_{0}-\mu}{\mu_{0}\cdot G}\\
 & =\phi\left(G\right)p_{0}\left(\frac{\mu-\mu_{0}}{\mu_{0}\cdot G}-\frac{\left(\ell+1\right)^{\bar{k}}\tau^{\nicefrac{1}{2}}}{k!\exp\left(nH\left(\hat{p},1-\hat{p}\right)\right)}\right).
\end{align*}
Remark that the first factor is $\phi\left(G\right)p_{0}$ is positive
and that 
\[
\frac{\left(\ell+1\right)^{\bar{k}}\tau^{\nicefrac{1}{2}}}{k!\exp\left(nH\left(\hat{p},1-\hat{p}\right)\right)}
\]
is a positive number that does not depend on $\mu_{0}.$ Further we
have
\[
\lim_{\mu_{0}\to0}\frac{\mu-\mu_{0}}{\mu_{0}\cdot G}=\infty\,\textrm{and}\,\lim_{\mu_{0}\to\infty}\frac{\mu-\mu_{0}}{\mu_{0}\cdot G}=0.
\]
Therefore it is sufficient to prove that $\frac{\mu-\mu_{0}}{\mu_{0}\cdot G}$
is a decreasing function of $\mu_{0}.$ 

We will use that $\frac{\mu-\mu_{0}}{\mu_{0}G}=\frac{p_{0}-\hat{p}}{\hat{p}\left(1-p_{0}\right)G},$
and prove that it is an increasing function of $p_{0}.$ The derivative
of $\frac{p_{0}-\hat{p}}{\left(1-p_{0}\right)G}$ can be calculated
as
\begin{align*}
\frac{\partial}{\partial p_{0}}\left(\frac{p_{0}-\hat{p}}{\left(1-p_{0}\right)G}\right) & =\frac{\left(1-p_{0}\right)G-\left(p_{0}-\hat{p}\right)\left(-G+\left(1-p_{0}\right)\frac{\partial G}{\partial p_{0}}\right)}{\left(1-p_{0}\right)^{2}G^{2}}\\
 & =\frac{1-\hat{p}-\left(p_{0}-\hat{p}\right)\frac{1-p_{0}}{2D}\frac{\partial D}{\partial p_{0}}}{\left(1-p_{0}\right)^{2}G}\\
 & =\frac{1-\hat{p}-\left(p_{0}-\hat{p}\right)\frac{1-p_{0}}{2D}\left(-\frac{\hat{p}}{p_{0}}+\frac{1-\hat{p}}{1-p_{0}}\right)\left(k+\ell\right)}{\left(1-p\right)^{2}G}\\
 & =\frac{1-\hat{p}-\frac{p_{0}-\hat{p}}{2D}\left(-\frac{\hat{p}\left(1-p_{0}\right)}{p_{0}}+1-\hat{p}\right)\left(k+\ell\right)}{\left(1-p_{0}\right)^{2}G}\\
 & =\frac{k+\ell}{\left(1-p_{0}\right)^{2}G^{3}}\left(\frac{2D}{k+\ell}\left(1-\hat{p}\right)-\frac{\left(p_{0}-\hat{p}\right)^{2}}{p_{0}}\right).
\end{align*}
The first factor changes sign from - to + at $p_{0}=\hat{p}$ and
the difference in the parenthesis equals 0 for $p_{0}=\hat{p}.$ The
derivative of the parenthesis is 
\begin{align*}
\frac{\partial}{\partial p_{0}}\left(\frac{2D}{k+\ell}\left(1-\hat{p}\right)-\frac{\left(p_{0}-\hat{p}\right)^{2}}{p_{0}}\right) & =2\frac{p_{0}-\hat{p}}{p_{0}\left(1-p_{0}\right)}\left(1-\hat{p}\right)-\frac{p_{0}2\left(p_{0}-\hat{p}\right)-\left(p_{0}-\hat{p}\right)^{2}}{p_{0}^{2}}\\
 & =\frac{\left(p_{0}-\hat{p}\right)^{2}\left(1+p_{0}\right)}{p_{0}^{2}\left(1-p_{0}\right)}\\
 & \geq0,
\end{align*}
so $\frac{D}{k+\ell}-\frac{\left(p_{0}-\hat{p}\right)^{2}}{2p_{0}\left(1-p_{0}\right)}$changes
sign from - to + at $p_{0}=\hat{p}.$

\section*{Proof of Theorem \ref{thm:InverseGaussian}}

A inverse Gaussian has density 
\[
f\left(x\right)=\left(\frac{\lambda}{\tau x^{3}}\right)^{\nicefrac{1}{2}}\exp\left(\frac{-\lambda\left(x-\mu\right)^{2}}{2\mu^{2}x}\right)
\]
with mean $\mu$ and shape parameter $\lambda.$ The variance function
is $V\left(\mu\right)=\mu^{3}/\lambda.$ The divergence of an inverse
Gaussian with mean $\mu_{1}$ from an inverse Gaussian with mean $\mu_{2}$
is
\begin{align*}
E\left[\ln\frac{\left(\frac{\lambda}{\tau x^{3}}\right)^{\nicefrac{1}{2}}\exp\left(\frac{-\lambda\left(x-\mu_{1}\right)^{2}}{2\mu_{1}^{2}x}\right)}{\left(\frac{\lambda}{\tau x^{3}}\right)^{\nicefrac{1}{2}}\exp\left(\frac{-\lambda\left(x-\mu_{2}\right)^{2}}{2\mu_{2}^{2}x}\right)}\right] & =E\left[\frac{\lambda\left(x-\mu_{2}\right)^{2}}{2\mu_{2}^{2}x}-\frac{\lambda\left(x-\mu_{1}\right)^{2}}{2\mu_{1}^{2}x}\right]\\
 & =\frac{\lambda}{2}E\left[\frac{x}{\mu_{2}^{2}}-\frac{2}{\mu_{2}}-\frac{x}{\mu_{1}^{2}}+\frac{2}{\mu_{1}}\right]\\
 & =\frac{\lambda}{2}\left(\frac{\mu_{1}}{\mu_{2}^{2}}-\frac{2}{\mu_{2}}-\frac{\mu_{1}}{\mu_{1}^{2}}+\frac{2}{\mu_{1}}\right)\\
 & =\frac{\lambda\left(\mu_{1}-\mu_{2}\right)^{2}}{2\mu_{1}\mu_{2}^{2}}.
\end{align*}
Therefore 
\[
g\left(x\right)=\frac{\lambda^{\nicefrac{1}{2}}\left(x-\mu\right)}{x^{\nicefrac{1}{2}}\mu}.
\]
Therefore the \emph{saddle-point approximation} is exact for the family
of inverse Gaussian distributions, i.e.
\[
f\left(x\right)=\frac{\phi\left(g\left(x\right)\right)}{\left(V\left(x\right)\right)^{\nicefrac{1}{2}}}.
\]

We have to prove that if $X$ is inverse Gaussian then $g\left(X\right)$
is majorize by a standard Gaussian. Equivalently we can prove that
if $Z$ is a standard Gaussian then $g^{-1}\left(Z\right)$ majorizes
the the inverse Gaussian. The density of $g^{-1}\left(Z\right)$ is
\begin{align*}
\tilde{f}\left(x\right) & =\phi\left(g\left(x\right)\right)\cdot g'\left(x\right)\\
 & =\phi\left(g\left(x\right)\right)\cdot\frac{x^{\nicefrac{1}{2}}\mu\lambda^{\nicefrac{1}{2}}-\nicefrac{1}{2}\cdot x^{-\nicefrac{1}{2}}\mu\lambda^{\nicefrac{1}{2}}\left(x-\mu\right)}{x\mu^{2}}\\
 & =\phi\left(g\left(x\right)\right)\cdot\lambda^{\nicefrac{1}{2}}\frac{x+\mu}{2x^{\nicefrac{3}{2}}\mu}\\
 & =\phi\left(g\left(x\right)\right)\cdot\frac{x+\mu}{2\left(V\left(x\right)\right)^{\nicefrac{1}{2}}\mu}.
\end{align*}

According to \cite[Lemma 5]{Harremoes2012} we can prove majorization
by proving that $\tilde{f}\left(x\right)/f\left(x\right)$is increasing.
We have 
\begin{align*}
\frac{\tilde{f}\left(x\right)}{f\left(x\right)} & =\frac{\phi\left(g\left(x\right)\right)\cdot\frac{x+\mu}{2\left(V\left(x\right)\right)^{\nicefrac{1}{2}}\mu}}{\frac{\phi\left(g\left(x\right)\right)}{\left(V\left(x\right)\right)^{\nicefrac{1}{2}}}}\\
 & =\frac{x+\mu}{2\mu},
\end{align*}
which is clearly an increasing function of $x$.

\section*{Poisson binomial distributions}

A \emph{Bernoulli sum} is a sum of independent Bernoulli random variables.
The distribution of a Bernoulli sum is called a \emph{Poisson binomial
distribution}. If the Bernoulli random variable $X_{i}$ has success
probability $p_{i}$ then the Bernoulli sum $S=\sum_{i=1}^{n}X_{i}$
has point probabilities
\[
P\left(S=j\right)=\prod_{i=1}^{n}\left(1-p_{i}\right)e_{j}\left(\frac{p_{1}}{1-p_{1}},\frac{p_{2}}{1-p_{2}},\dots,\frac{p_{n}}{1-p_{n}}\right)
\]
where $e_{j}$ denote the symmetric polynomials
\begin{align*}
e_{0}\left(x_{1},x_{2},\dots,x_{n}\right) & =1\\
e_{1}\left(x_{1},x_{2},\dots,x_{n}\right) & =x_{1}+x_{2}+\dots+x_{n}\\
e_{2}\left(x_{1},x_{2},\dots,x_{n}\right) & =x_{1}x_{2}+x_{1}x_{3}+\dots x_{n-1}x_{n}\\
 & \vdots\\
e_{n}\left(x_{1},x_{2},\dots,x_{n}\right) & =x_{1}x_{2}\cdot\dots\cdot x_{n}
\end{align*}
The probability generating function can be written as
\begin{align*}
M\left(t\right) & =\sum_{j=0}^{n}P\left(S=j\right)t^{j}\\
 & =\sum_{j=0}^{n}\prod_{i=1}^{n}\left(1-p_{i}\right)e_{j}\left(\frac{p_{1}}{1-p_{1}},\frac{p_{2}}{1-p_{2}},\dots,\frac{p_{n}}{1-p_{n}}\right)t^{j}\\
 & =\prod_{i=1}^{n}\left(1-p_{i}\right)\sum_{j=0}^{n}e_{j}\left(\frac{p_{1}}{1-p_{1}},\frac{p_{2}}{1-p_{2}},\dots,\frac{p_{n}}{1-p_{n}}\right)t^{j}\\
 & =\prod_{i=1}^{n}\left(1-p_{i}\right)\prod_{i=1}^{n}\left(\frac{p_{i}}{1-p_{i}}+t\right).
\end{align*}
We see that $t$ is a root of the probability generating function
if and only if $t=-\frac{p_{i}}{1-p_{i}}$ for some value of $i$.
Hence the probability generating function of a Bernoulli sum has $n$
real roots counted with multiplicity. All roots of a probability generating
function are necessarily non-positive so if a probability generating
function has $n$ real roots then they are automatically non-positive
and the function must be the probability generating function for a
Bernoulli sum with success probabilities that can be calculated as
$p_{i}=-\frac{t}{1-t}$. Using this type of properties it has been
shown \cite{Barbour1992} that all hypergeometric distributions are
Poisson binomial distributions.
\begin{prop}
If $P$ is a Poisson binomial distribution then all distributions
in the exponential family generated by $P$ are Poisson binomial distributions.\end{prop}
\begin{IEEEproof}
Assume that the random variable $S$ has distribution $P$. Then the
probability generating function of $P_{\beta}$ is
\begin{align*}
\sum_{j=0}^{n}P_{\beta}\left(S=j\right)t^{j} & =\sum_{j=0}^{n}\frac{\exp\left(\beta j\right)}{Z\left(\beta\right)}P\left(S=j\right)t^{j}\\
 & =\frac{1}{Z\left(\beta\right)}\sum_{j=0}^{n}P\left(S=j\right)\left(\exp\left(\beta\right)t\right)^{j}
\end{align*}
so $t$ is root of the probability generating function of $S$ if
and only $\exp\left(\beta\right)t$ is root of the probability generating
function of $P_{\beta}.$
\end{IEEEproof}
The new success probabilities can be calculated from the old ones
by
\begin{align*}
\tilde{p}_{i} & =-\frac{\exp\left(\beta\right)\left(-\frac{p_{i}}{1-p_{i}}\right)}{1-\exp\left(\beta\right)\left(-\frac{p_{i}}{1-p_{i}}\right)}\\
 & =\frac{\exp\left(\beta\right)p_{i}}{1+\left(\exp\left(\beta\right)-1\right)p_{i}}.
\end{align*}

The following result is well-known, but we give the proof for completeness.
\begin{thm}
\label{thm:median}Let $S$ denote a Bernoulli sum with $E\left[S\right]=m\in\mathbb{N}_{0}$.
Then the median is $m$, i.e. $P\left(S<m\right)<\nicefrac{1}{2}<P\left(X\leq m\right)$. \end{thm}
\begin{IEEEproof}
By symmetry it is sufficient to prove that $P\left(S<m\right)<\nicefrac{1}{2}.$
Let $S=X_{1}+X_{2}+\tilde{S}.$ Then 
\begin{align*}
P\left(S<m\right) & =P\left(\tilde{S}<m-2\right)+P\left(\tilde{S}=m-2\right)P\left(X_{1}+X_{2}\leq1\right)+P\left(\tilde{S}=m-1\right)P\left(X_{1}+X_{2}=0\right)\\
 & =P\left(\tilde{S}<m-2\right)+P\left(\tilde{S}=m-2\right)\left(\left(1-p_{1}\right)\left(1-p_{2}\right)+\left(p_{1}\left(1-p_{2}\right)+\left(1-p_{1}\right)p_{2}\right)\right)\\
 & \,\,+P\left(\tilde{S}=m-1\right)\left(1-p_{1}\right)\left(1-p_{2}\right)\\
 & =P\left(\tilde{S}<m-2\right)+P\left(\tilde{S}=m-2\right)\left(1-p_{1}p_{2}\right)+P\left(\tilde{S}=m-1\right)\left(1-\left(p_{1}+p_{2}\right)+p_{1}p_{2}\right)\\
 & =P\left(\tilde{S}<m\right)+\left(P\left(\tilde{S}=m-1\right)-P\left(\tilde{S}=m-2\right)\right)p_{1}p_{2}-P\left(\tilde{S}=m-1\right)\left(p_{1}+p_{2}\right)
\end{align*}
If we fix $p_{1}+p_{2}$then we see that $P\left(S<m\right)$ is maximal
when $p_{1}=p_{2}$ because $P\left(\tilde{S}=m-1\right)\geq P\left(\tilde{S}=m-2\right).$
Using this type of argument we see that $P\left(S<m\right)$ is maximal
when all success probabilities are equal, i.e. the distribution of
$S$ is binomial. 
\end{IEEEproof}

\section*{Numerical calculations}

In this section all computations are performed in in R programming
language. First we will prove the conjecture in some special cases.
These special cases are of independent interest but by avoiding them
we are also able to simplify the program for numerical checking of
the conjecture and these simplifications will decrease the execution
time of the program. First we prove the conjecture i the case when
one of the cells are empty. Without loss of generality we may assume
that $x=0.$
\begin{thm}
For a $2\times2$ contingency table let $G\left(X\right)$ denote
the signed log-likelihood. Then $\Pr\left(X<0\right)<\Pr\left(Z\leq G\left(0\right)\right)<\Pr\left(X\leq0\right)$
where $Z$ is a standard Gaussian.\end{thm}
\begin{IEEEproof}
The left part of the double inequality is trivially fulfilled so we
just have to prove that 
\[
\Pr\left(Z\leq G\left(0\right)\right)<\Pr\left(X=0\right).
\]
A large deviation bound of the standard Gaussian leads to
\[
\Pr\left(Z\leq G\left(0\right)\right)\leq\exp\left(-N\cdot I\right)
\]
and 
\[
\Pr\left(X=0\right)=\frac{\left(N-n\right)!\left(N-r\right)!}{\left(N-n-r\right)!N!}
\]
and $N\cdot I=\left(N-n\right)\ln\left(N-n\right)+\left(N-r\right)\ln\left(N-r\right)-\left(N-n-r\right)\ln\left(N-n-r\right)-N\ln\left(N\right)$.
Hence we have to prove that
\[
\frac{\left(N-n\right)^{N-n}\left(N-r\right)^{N-r}}{\left(N-n-r\right)^{N-n-r}N^{N}}\leq\frac{\left(N-n\right)!\left(N-r\right)!}{\left(N-n-r\right)!N!}
\]
or equivalently
\[
\frac{\left(N-n-r\right)!N!}{\left(N-n-r\right)^{N-n-r}N^{N}}\leq\frac{\left(N-n\right)!\left(N-r\right)!}{\left(N-n\right)^{N-n}\left(N-r\right)^{N-r}}.
\]
To prove this define $k=n+r$ so that we have to prove that
\[
\frac{\left(N-k\right)!N!}{\left(N-k\right)^{N-k}N^{N}}\leq\frac{\left(N-n\right)!}{\left(N-n\right)^{N-n}}\cdot\frac{\left(N-k+n\right)!}{\left(N-k+n\right)^{N-n}}.
\]
For $n=0$and $n=k$ this is obvious so we just have to prove that
the function 
\[
f\left(n\right)=\frac{\left(N-n\right)!}{\left(N-n\right)^{N-n}}\cdot\frac{\left(N-k+n\right)!}{\left(N-k+n\right)^{N-n}}
\]
is first increasing and then decreasing. We have 
\begin{align*}
\frac{f\left(n+1\right)}{f\left(n\right)} & =\frac{\frac{\left(N-n-1\right)!}{\left(N-n-1\right)^{N-n-1}}\cdot\frac{\left(N-k+n+1\right)!}{\left(N-k+n+1\right)^{N-n+1}}}{\frac{\left(N-n\right)!}{\left(N-n\right)^{N-n}}\cdot\frac{\left(N-k+n\right)!}{\left(N-k+n\right)^{N-n}}}\\
 & =\frac{\left(1+\frac{1}{N-n-1}\right)^{N-n-1}}{\left(1+\frac{1}{N-k+n}\right)^{N-k+n}}.The
\end{align*}
result follows because $\left(1+\frac{1}{x}\right)^{x}$ is a decreasing
function.\end{IEEEproof}
\begin{thm}
The intersection conjecture for hypergeometricd distributions is true
if $x=nr/N.$\end{thm}
\begin{IEEEproof}
In this case $G\left(x\right)=0$and $\Phi\left(G\left(x\right)\right)=\nicefrac{1}{2}.$
In this case the intersection conjecture states that
\[
\Pr\left(X<x\right)\leq\nicefrac{1}{2}\leq\Pr\left(X\leq x\right)
\]
i.e. $x$ is is the median of the distribution but according to Theorem
\ref{thm:median} this holds for any Poisson binomial distribution.
\end{IEEEproof}
Due to symmetry of the problem we only have to check the conjecture
for $x<n\cdot r/N.$ The following program tests that the inequality
is correct for any contingency table with a sample size less than
or equal to 200. The \emph{for loop} is slow in the R program so if
the computer is not vary fast one may want to chop the outer loop
into pieces (say 4 to 50, 51 to 100 etc.)

\begin{lstlisting}
xlnx <- function(x){x*log(x)}
mutual <- function(x,n,k,r){xlnx(x)+xlnx(n-x)+xlnx(r-x)+xlnx(k-r+x)
                     -xlnx(n)-xlnx(k)-xlnx(r)-xlnx(n+k-r)+xlnx(n+k)} 
	# calculates the mutual information
g <- function(x,n,k,r){-sqrt(2*mutual(x,n,k,r))} 
	# calculates the signed log-likelihood 
f1 <- function(x,n,k,r){phyper(x-1,n,k,r) < pnorm(g(x,n,k,r)) 
                       & pnorm(g(x,n,k,r)) < phyper(x,n,k,r) } 
	# checks if the conjecute is satisfied for specific values
f2 <- function(n,k,r){x=max(1,r-k+1); while (x < n*r/(n+k)) 
                            {w <- w&f1(x,n,k,r);x<-(x+1)};w}  
	# loop that checks the conjecture for all values 
	# of x less than the mean 
f3 <- function(n,k){for (r in 2:(n+k-2)) w <- w&f2(n,k,r);w}  
	# loop that checks the conjecture for all values of the row sum
f4 <- function(tot){for (n in 2:(tot-2)) w <- w&f3(n,tot-n);w}   
	# loop that checks the conjecture for all values of the column sum
w=TRUE	
	# set the initial value of the variable w to true
for (tot in 4:200) w <- w&f4(tot); w 	
	# run throug all sample sizes from 4 to 200
\end{lstlisting}

The result is TRUE.

Next we want to check the intersection conjecture in the cases where
one of the bins contain less than 10 observations. Without loss of
generality we may assume that $x\leq9$. If $n\leq0.1\cdot N$ or
if $r\leq0.1\cdot N$ then we can make a binomial approximation. Assume
that $\lambda=n\cdot r/N\geq10.$ Then $\Pr\left(X\leq x\right)\leq Po\left(\lambda,x\right)\leq Po\left(\lambda,9\right).$
Now $Po\left(22.7,9\right)<0.001$ so if $\lambda\geq22.7$ and $x\leq9$
then the null hypothesis of independence can be rejected using this
upper bound by a Poisson distribution. Therefore it is only necessary
to check the intersection inequality for $\lambda\leq22.7$. Since
$\lambda=n\cdot\left(r/N\right)\geq0.1\cdot n$we have $n\leq227,$
and similarly $r\leq227.$ Finally we have $\lambda=\nicefrac{n}{N}\cdot\nicefrac{r}{N}\cdot N\geq0.01\cdot N$
implying that $N\leq2270.$ Hence we will check all cases where $1\leq x\leq\min\left\{ 9,n\cdot r/N\right\} $
and $2\leq n\leq227$ and $2\leq r\leq227$ and $4\leq N\leq2270.$

The following program tests that the conjecture is correct.

\begin{lstlisting}
xlnx <- function(x){x*log(x)}
mutual <- function(x,n,k,r){xlnx(x)+xlnx(n-x)+xlnx(r-x)+xlnx(k-r+x)
                     -xlnx(n)-xlnx(k)-xlnx(r)-xlnx(n+k-r)+xlnx(n+k)}
g <- function(x,n,k,r){-sqrt(2*mutual(x,n,k,r))}
f1 <- function(x,n,k,r){pnorm(g(x,n,k,r)) < phyper(x,n,k,r)
                        & phyper(x-1,n,k,r) < pnorm(g(x,n,k,r))}
f2 <- function(n,k,r){x=max(1,r-k+1); while (x<10 && x < n*r/(n+k)) 
                                    {w <- w&f1(x,n,k,r);x<-(x+1)};w}  
f3 <- function(n,k){for (r in 2:(min(n+k-2,227,22.7*(1+k/n)))) 
                                w <- w&f2(n,k,r);w}  
f4 <- function(tot){for (n in 2:(min(tot-2,227))) 
                             w <- w&f3(n,tot-n);w} 
w=TRUE
for (tot in 4:2270) w <- w&f4(tot); w
\end{lstlisting}

The result is TRUE.

\section*{QQ-plots for other contingency tables}

As we have seen that the distribution of the signed log-likelihood
is close to Gaussian for contingency tables with one degree of freedom.
If we square the signed log-likelihood we get that mutual information
is close to a $\chi^{2}$-distribution but if the contingency table
is not symmetric the steps in the negative and positive direction
may lead to a confusing interference pattern with fluctuations away
from a straight line as illustrated in the next example.
\begin{example}
Consider a contingency table with expected counts as in table \ref{tab:ContingencyTableDF=00003D2-1}.
\begin{table}
\begin{centering}
\begin{tabular}{|c|c|c|}
\hline 
4 & 16 & 20\tabularnewline
\hline 
16 & 64 & 80\tabularnewline
\hline 
20 & 80 & 100\tabularnewline
\hline 
\end{tabular}
\par\end{centering}

\caption{\label{tab:ContingencyTableDF=00003D2-1}Contingency table with expected
values and marginal counts.}
\end{table}
In Figure \ref{fig:QQasym} we see that there are steps of different
sizes, where some come from positive values of the signed log-likelihood
and some come from negative values.
\begin{figure}
\centering{}\includegraphics[scale=0.4]{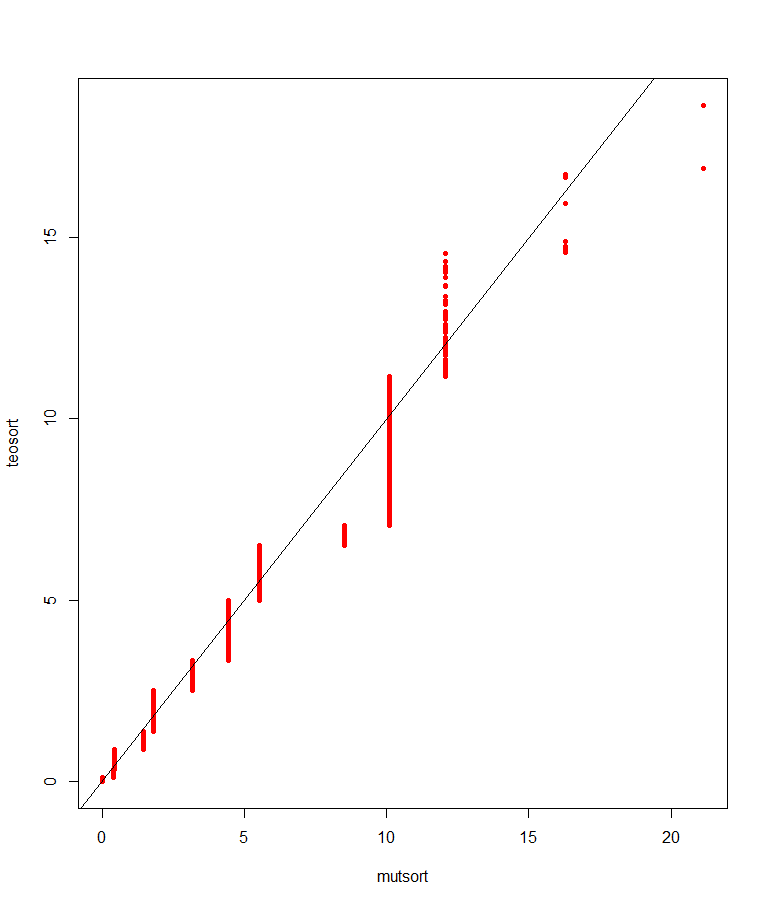}\caption{\label{fig:QQasym}A QQ-plot of a sample of the $G^{2}$-statistic
from the contingency table \ref{tab:ContingencyTableDF=00003D2}versus
a sample from a $\chi^{2}$-distribution with two degrees of freedom.
the sample size in this simulation is 100 000.}
\end{figure}

To deminish such fluctuations we can make an asymmetric continuity
correction as proposed by Gy{\"e}rfi, Harremo{\"e}s and Tusn{\'a}dy \cite{Gyorfi2012}.
Such continuity correction will lead to a much closer fit to a $\chi^{2}$-distribution,
but it implies that the plugin estimate of mutual information has
to be replaced by another expression and is thus not within the scope
of this paper. Therefore the following examples will focus on symmetric
contingency tables.
\end{example}
One may wonder if our observations are also true if there is more
than one degree of freedom. It turns out out to be a more complicated
question with no definite answer. We will illustrate this by two examples.
\begin{example}
In this example we consider the contingency table \ref{tab:ContingencyTableDF=00003D2}.
\begin{table}
\begin{centering}
\begin{tabular}{|c|c|c|c|}
\hline 
18 & 18 & 18 & 54\tabularnewline
\hline 
18 & 18 & 18 & 54\tabularnewline
\hline 
36 & 36 & 36 & 108\tabularnewline
\hline 
\end{tabular}
\par\end{centering}

\caption{\label{tab:ContingencyTableDF=00003D2}Contingency table with expected
values and marginal counts.}
\end{table}

\end{example}
Due to the symmetry we do not get signicant interference between fluctuations
in different direction. Both the distribution of the $\chi^{2}$-statistic,
the distribution of mutual information, and the theoretical $\chi^{2}$-distribution
have been simulated with a sample size equal to one million. In Figure
\ref{fig:A-QQ-plotMutvsTeo} we see that the intersection property
holds for mutual information , and in Figure \ref{fig:A-QQ-plotchivsTeo}that
the $\chi^{2}$-statistic systematically gives values that are too
small for large deviation levels.

\begin{figure}
\begin{centering}
\includegraphics[scale=0.4]{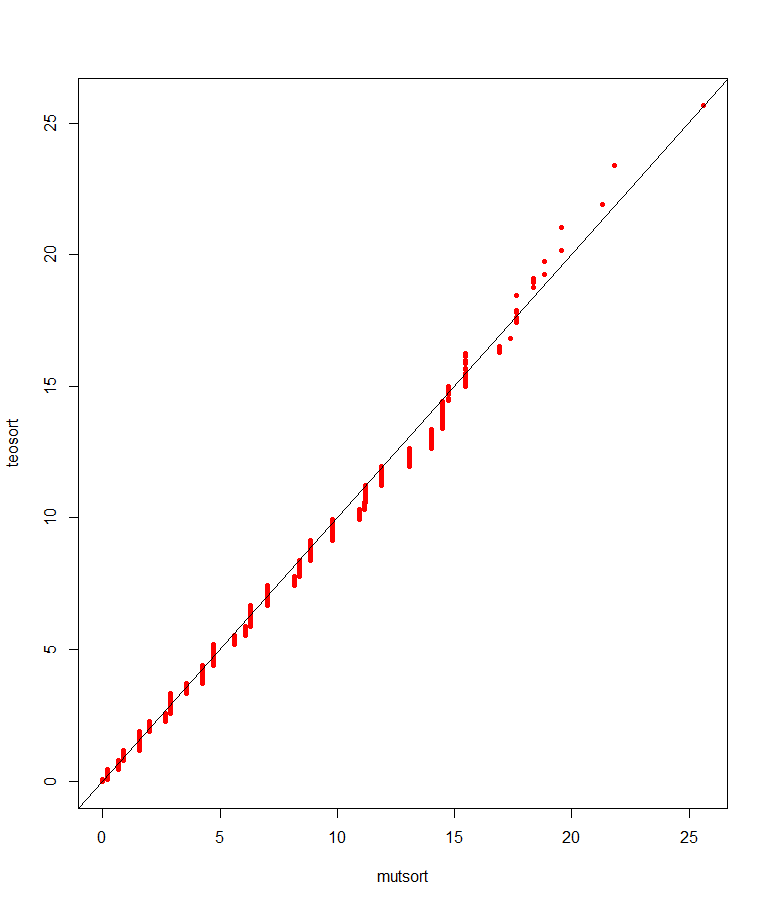}
\par\end{centering}

\caption{\label{fig:A-QQ-plotMutvsTeo}A QQ-plot of a sample of the $G^{2}$-statistic
from the contingency table \ref{tab:ContingencyTableDF=00003D2}versus
a sample from a $\chi^{2}$-distribution with two degrees of freedom.
the sample size in this simulation is 100 000.}

\end{figure}

\begin{figure}
\begin{centering}
\includegraphics[scale=0.4]{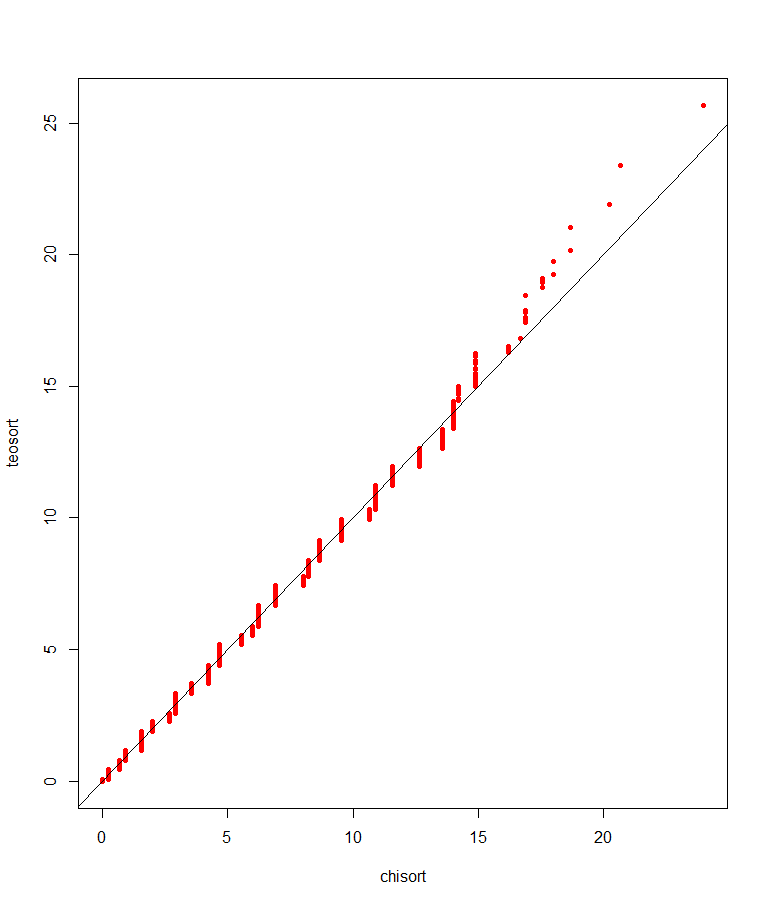}
\par\end{centering}

\caption{\label{fig:A-QQ-plotchivsTeo}A QQ-plot of a sample of the $\chi^{2}$-statistic
from the contingency table \ref{tab:ContingencyTableDF=00003D2}versus
a sample from a $\chi^{2}$-distribution with two degrees of freedom.
the sample size in this simulation is 100 000.}

\end{figure}

\begin{example}
For the contingency table \ref{tab:ContMangefrih} the sample size
is relatively small compared with the number of degrees of freedom.
In Figure \ref{fig:A-QQ-plotMagefrih} we see that the $\chi^{2}$-statistic
has a sligthly smaller vaue than the theoretical distribution while
the mutual information is systematically too large for large deviations
in Figure \ref{fig:A-QQ-plotchiMange}.
\begin{figure}
\begin{centering}
\includegraphics[scale=0.4]{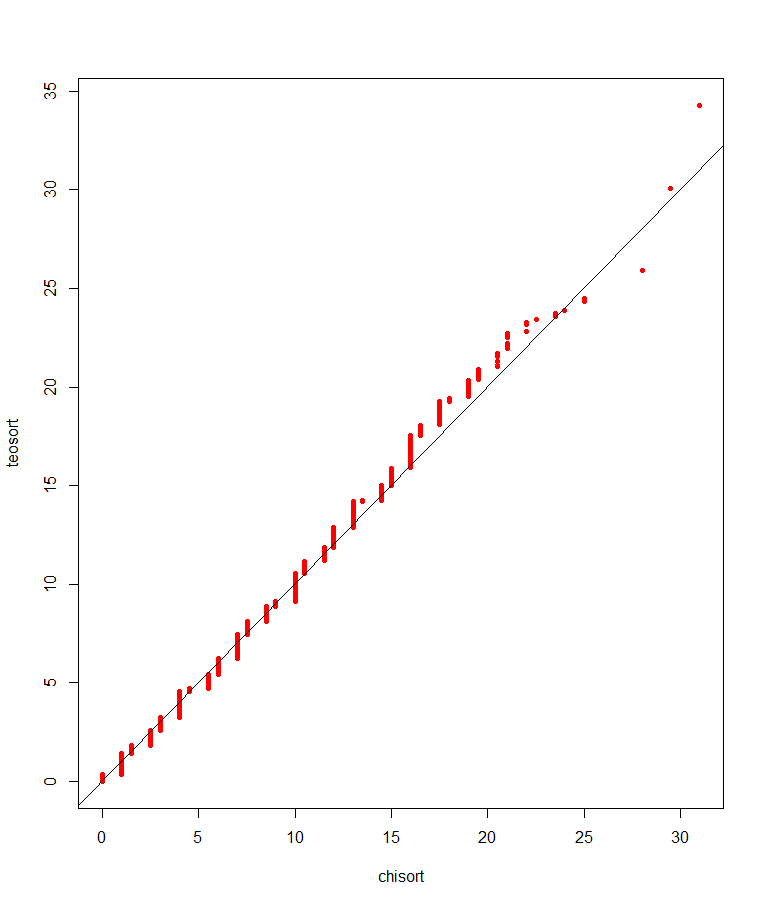}
\par\end{centering}

\caption{\label{fig:A-QQ-plotMagefrih}A QQ-plot of a sample of the $\chi^{2}$-statistic
from the contingency table \ref{tab:ContMangefrih} versus a sample
from a $\chi^{2}$-distribution with two degrees of freedom. the sample
size in this simulation is 100 000.}

\end{figure}

\end{example}
\begin{center}
\begin{table}
\begin{centering}
\begin{tabular}{|c|c|c|c|}
\hline 
4 & 4 & 4 & 12\tabularnewline
\hline 
4 & 4 & 4 & 12\tabularnewline
\hline 
4 & 4 & 4 & 12\tabularnewline
\hline 
12 & 12 & 12 & 36\tabularnewline
\hline 
\end{tabular}
\par\end{centering}

\caption{\label{tab:ContMangefrih}Contingency table with expected values and
marginal counts.}
\end{table}

\par\end{center}

\begin{figure}
\begin{centering}
\includegraphics[scale=0.4]{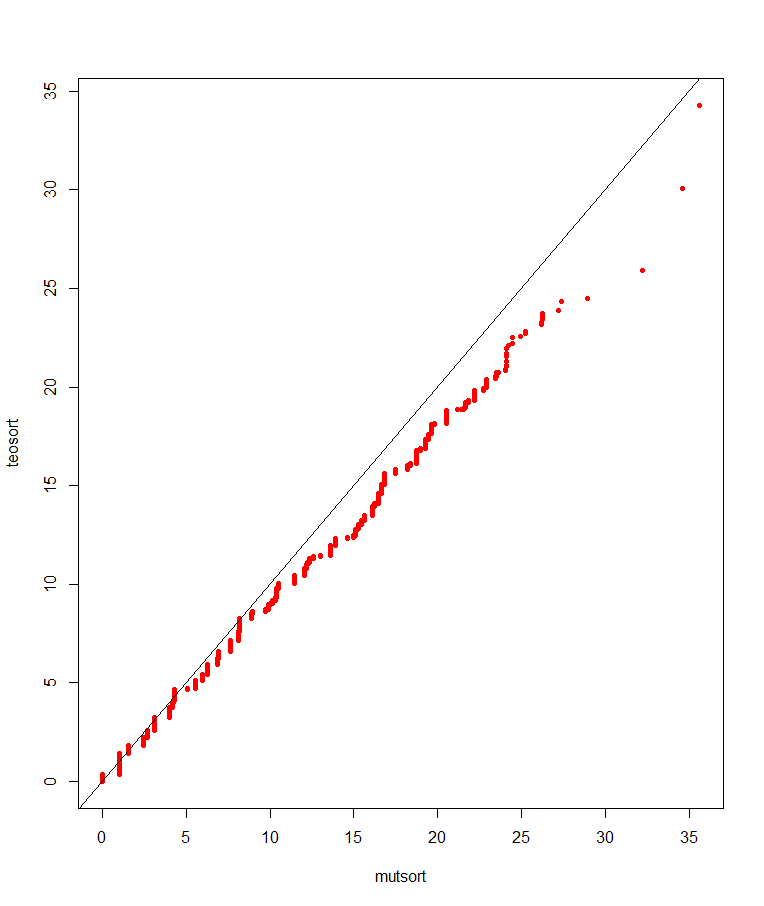}\caption{\label{fig:A-QQ-plotchiMange}A QQ-plot of a sample of the $\chi^{2}$-statistic
from the contingency table \ref{tab:ContMangefrih}versus a sample
from a $\chi^{2}$-distribution with two degrees of freedom. the sample
size in this simulation is 100 000.}

\par\end{centering}

\end{figure}

The tendency is that mutual information is better for large deviations
and $\chi^{2}$ is better if the sample size is large compared with
the degrees of freedom. The last problem can be corrected by 'smoothing'
by a rate distortion test as suggested in \cite{Harremoes2008g},
but both continuity corrections and rate distortion tests just tell
that there are good alternative to $\chi^{2}$ and mutual information.
Looking at the behavior of such alternatives is not relevant if the
question how to choose between $\chi^{2}$and mutual information.
We see that if there are more than one degree of freedom several effects
point in opposite directions, so it might be more important to reash
a better understanding of the situation where there is one degree
of freedom before trying to generalize the results about one degree
of freedom.
\end{document}